\documentclass[12pt]{article}
\usepackage{epsfig,amsfonts,amsmath,amssymb,amscd,authblk}

\newcommand{\bfe}{\mbox{\boldmath $e$}}

\newcommand{\bfu}{\mbox{\boldmath $u$}}

\newcommand\Tstrut{\rule{0pt}{2.6ex}}         
\newcommand\Bstrut{\rule[-0.9ex]{0pt}{0pt}}  
\usepackage[table]{xcolor}
\usepackage{subcaption}
\usepackage{graphicx} 
\usepackage{setspace}
\usepackage{hyperref}

\begin{document}
	
	\title {Axisymmetric deformation of compressible, nearly incompressible, and incompressible  thin layers between two rigid surfaces\footnote{Published in {\it International Journal of Solids and Structures}, vol. 214–215, pp. 61-73, 2021. https://doi.org/10.1016/j.ijsolstr.2020.12.002}}
	
	\author[a]{Alexander B. Movchan}
	\author[b]{Kirill R. Rebrov}
	\author[b,c]{\\Gregory J. Rodin\thanks{For correspondence contact at gjr@oden.utexas.edu}
	}
	\bigskip
	
	\affil[a]{Department of Mathematical Sciences\\ University of Liverpool, L69 7ZL
		UK}
	\affil[b]{Oden Institute for Computational Engineering and Sciences, The University of Texas at Austin\\ Austin, TX 78712 USA}
		\affil[c]{Department of Aerospace Engineering and Engineering Mechanics,
		The University of Texas at Austin\\ Austin, TX 78712 USA } 

	\date{}
	
	\maketitle
	\begin{abstract} \noindent
		Accurate asymptotic solutions are presented for axisymmetric deformation of thin layers constrained by either two rigid plates or two rigid spheres. Those solutions are developed using Saint-Venant's principle and the layer thinness as the only assumptions. The solutions are valid in the entire range of Poisson's ratios, and allow one to distinguish among compressible,  nearly incompressible, and  incompressible layers. That classification involves both material and geometric parameters. 
	\end{abstract}
	\bigskip

	\pagebreak
	
\section{Introduction}
This paper is concerned with analysis of axisymmetric deformation of thin layers between two rigid surfaces. This problem has been studied in depth in the fluid mechanics literature for flat surfaces (plates), because  of its significance to rheometry; see \cite{Engmann++} for references which include a wide range of constitutive models and conditions along the layer-plate interfaces. According to  \cite{Engmann++}, the first asymptotic solution, which takes advantage of the layer thinness, was published by Jo\v{z}ef Stefan in 1874.

\begin{figure}[h]
	\begin{center} 
		\includegraphics[width=\textwidth]{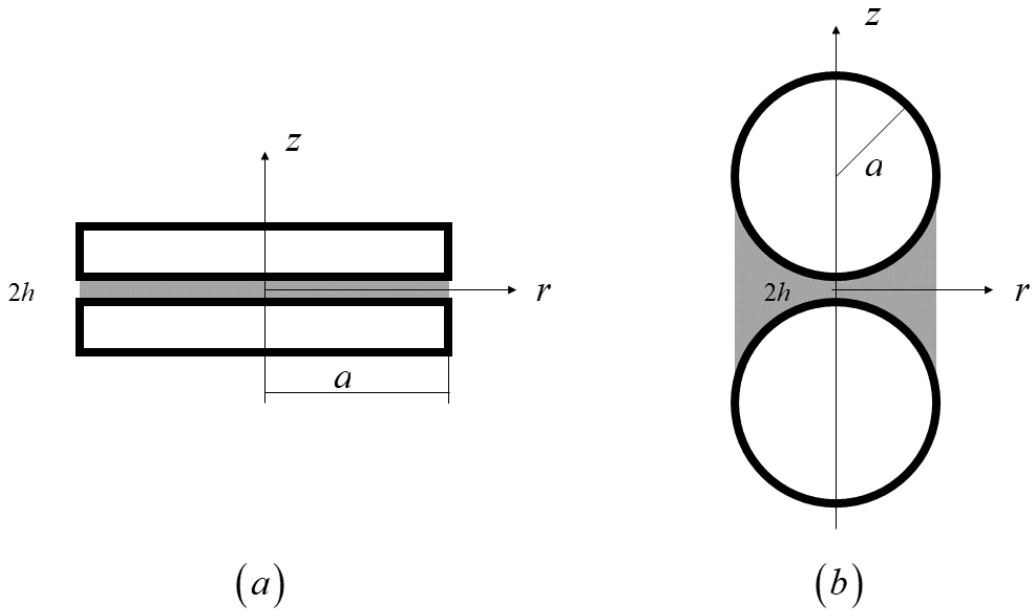}\\
		\caption{Thin layers between two rigid surfaces: (a) plates and (b) spheres.}		\label{Fig1} 
	\end{center}
\end{figure}

\bigskip
\noindent
In the basic setting, when perfect bonding and incompressibility are assumed, squeezing of a thin Newtonian fluid layer by two plates is described by the velocities
\begin{equation}\label{fluids radial velocity}
v_r (r,z) = \frac{3 r V
	\left(h^2-z^2\right)}{4
	h^3}\:,
\end{equation} \smallskip
\begin{equation}\label{fluids axial velocity}
v_z (r,z) = \frac{V z \left(z^2-3
	h^2\right)}{2 h^3}
\:,
\end{equation} \\
and pressure 
\begin{equation}\label{fluids pressure}
p(r,z) = \frac{\mu  V \left(3 a^2+2
	h^2-3 r^2+6 z^2\right)}{4
	h^3}
\:.
\end{equation} \\
In these equations, $a$ and $h$ are the layer radius and half-thickness, respectively, and $\mu$ is the shear viscosity of the fluid. The cylindrical coordinates are naturally aligned with the layer ($0\leq r < a$ and $-h<z<h$), and each plate moves toward the other with velocity $ V$ (Fig.~\ref{Fig1}a).

\bigskip
\noindent
The fields in (\ref{fluids radial velocity}--\ref{fluids pressure}) do not result in traction-free boundary conditions on the cylindrical surface, as 
\begin{equation*}\label{fluid strong BC r}
\sigma_{rr}(a,z)=-p(a,z)+2\mu \frac{\partial v_r (r,z)}{\partial r}_{|r=a}=\frac{\mu  V \left(h^2-3
	z^2\right)}{h^3}
\end{equation*}
and
\begin{equation*}\label{fluid strong BC z}
\sigma_{rz}(a,z)=\mu \left[ \frac{\partial v_r (r,z)}{\partial z} +  \frac{\partial v_z (r,z)}{\partial r}\right]_{|r=a} =-\frac{3 a \mu  V z}{2 h^3}\:.
\end{equation*} \\
Nevertheless, these boundary conditions are satisfied  
in Saint-Venant's sense, as both resultants are equal to zero:
$$
\int_{-h}^{h}\sigma_{rr}(a,z) \: {\rm d}z=\int_{-h}^{h}\sigma_{rz}(a,z) \: {\rm d}z=0\:.
$$
This means that the fields in (\ref{fluids radial velocity}--\ref{fluids pressure}) are valid only at sufficiently large distances away from the cylindrical surface, which implies that (\ref{fluids radial velocity}--\ref{fluids pressure}) are meaningful only when $a\gg h$, that is, the layer must  be thin. This Saint-Venant's setting is central to our and all preceding developments.

	\bigskip
	\noindent
	 In the mathematical literature, Saint-Venant's principle is associated with asymptotic solutions of boundary-value problems defined on unbounded domains. In particular, for infinite strips, Saint-Venant's principle is directly connected to analysis of boundary layers \cite{oleinik}. In this work, we exploit this connection in the context of method of compound asymptotic approximations \cite{Maz'ya, Movchans1}, particularly effective for thin domains. This method implies that rigid confining surfaces, which subject thin layers to Dirichlet boundary conditions, induce boundary layers characterized by exponential decays away from the cylindrical surface. That is, Saint-Venant's principle is fully expected to hold for the problems of interest.

\bigskip
\noindent
 Let us emphasize two remarkable properties of the fields in (\ref{fluids radial velocity}--\ref{fluids pressure}):
\begin{itemize}
	\item 
The pressure at the center is
$$
p(0,0)\approx \frac{3\mu  V }{h} \times \left( \frac{a}{2h} \right)^2  \:.
$$
The first fraction in this expression would be the pressure if the fluid were allowed to slip freely along the plates. Thus the no-slip condition on the fluid-plate interfaces results in a dramatic pressure build-up near the center.
\item 
The pressure on the plates at $r=0$  exceeds that at the cylindrical surface  by two orders of magnitude, and the maximum deviatoric stress by one order of magnitude. As a result, one can accurately calculate the forces acting on the plates based on (\ref{fluids pressure}) alone, without taking into account the deviatoric stresses. 
\end{itemize}
These properties of the pressure field and the simplicity of  (\ref{fluids radial velocity}--\ref{fluids pressure}) has resulted in numerous approximate solutions involving non-Newtonian fluids and partial-slip boundary conditions along the plate-fluid interfaces \cite{Engmann++}.

\bigskip
\noindent
In this paper, we are primarily interested in solid rather than fluid layers. In the simplest setting, when the layer is linear elastic, and the bonding is perfect, the key difference between solids and fluids, at least as far as mathematics is concerned, is compressibility. This difference is so significant that the basic asymptotic ansatz ubiquitous to analysis of thin  incompressible layers must be re-examined; this issue will be addressed in Section 2. This may explain why existing analyses of compressible layers are based on various additional assumptions, in the spirit of Bernoulli-Euler beam theory or Reynold's lubrication theory. Earlier solutions following this path are reviewed in \cite{Pinarbasi2008}. There it is observed that most of those approximate solutions are based on assumptions inspired by (\ref{fluids radial velocity}--\ref{fluids pressure}), in part because originally the problem was motivated by studies of rubber, whose elastic response is very close to being incompressible. We will discuss those solutions in Section 2, once a proper mathematical setting has been introduced. More recent studies \cite{H2002++, QL2015++, RAS2} adopt less restrictive assumptions, which may be useful for not very thin layers. Of course, the problem does not pose significant challenges for finite element computations, as long as the layer is not too close to being incompressible \cite{QL2015++}. 

\bigskip
\noindent
The significant role of compressibility can be further explored by considering a closely-related problem involving a thin layer between two rigid spheres. For this problem,  asymptotic solutions are available for both compressible \cite{PTK} and incompressible \cite{Jeffrey} cases. Equations for the  squeezing force derived by those authors expose a chasm. In particular, for  compressible layers, 
$$
F^c \sim \frac{1}{1-2\nu}\log\frac{a}{h}\:,
$$
where $a$ is the sphere radius, $2h$ is the smallest layer thickness, and $\nu$ is Poisson's ratio (Fig.~\ref{Fig1}b). This equation  becomes problematic as $\nu \rightarrow 1/2$, that is, in the limit as the layer material approaches incompressibility.  For incompressible layers, 
$$
F^i \sim \frac{a}{h}\:.
$$
Thus it is unclear which of the two equations to use  for nearly incompressible solids like rubber, with  Poisson's ratio close to one half.   Another interesting feature of the two solutions is that for compressible layers 
$$
\frac{u_r}{u_z}={\cal O}\left(\frac{h}{a}\right)\:,
$$
whereas for incompressible layers 
$$
\frac{u_r}{u_z}={\cal O}\left(\frac{a}{h}\right)\:.
$$
\noindent
The chasm between the two asymptotic solutions is easy to explain, but not resolve, in terms of asymptotic analysis. For nearly incompressible layers, the problem involves two small parameters, $h/a$ and  ${1-2\nu}$, and a proper approach must take this into account. In contrast, the analyses of \cite{PTK} and \cite{Jeffrey} consider only the extremes, each involving only one small parameter  $h/a$. In this regard, it is not surprising that the two solutions do not match. Also, let us mention that, besides similarities, the problems for flat and spherical constraining surfaces are rather different. In particular, for flat surfaces, the curvature $\kappa=0$, whereas for spherical surfaces $\kappa a =1$. Thus the two problems can be regarded as two extreme cases of a class of problems where one may have two rather than one dimensionless geometric parameters, representing the layer thickness, constraint size, and curvature. Additional differences will be exposed as the solutions are being developed. 

\bigskip
\noindent
The goal of this paper is to construct new approximate solutions   for axisymmetric  stretching of thin linear elastic layers constrained by two rigid plates or two rigid spheres in the entire range of Poisson's ratio  $-1< \nu \leq 1/2$. These solutions are based on two assumptions: ($i$) the layer thinness and ($ii$) Saint-Venant's form of the boundary conditions on the cylindrical surface. In contrast, existing approximate solutions rely on these two assumptions plus some other assumptions. Note that the thinness assumption is not explicitly stated in \cite{QL2015++, RAS2}, but it is implied once Saint-Venant's form of the boundary conditions on the cylindrical surface is adopted. We will demonstrate that our solutions are accurate in the entire range of Poisson's ratio  $-1< \nu \leq 1/2$. Further, we will exploit the two solutions to properly define compressible, intermediate, and incompressible responses, which take  into account not just the proximity of Poisson's ratio to one half,  but also the geometry.

\bigskip
\noindent
The remainder  of this paper consists of five sections. In Sections 2 and 3, we present  complete leading order asymptotic solutions for stretching of linear elastic layers constrained by two rigid plates (Section 2), and by two rigid spheres (Section 3). In Section 4, we compare existing single-parameter and new two-parameter asymptotic solutions, and define compressible, intermediate, and incompressible regimes. In Section 5, we develop an alternative approach, based on asymptotic series, which provides a better understanding of transitions from compressible to  intermediate, and from intermediate to incompressible regimes. In Section 6, we discuss various connections of our work with related problems for thin layers, and outline
possible extensions.

\section{Thin layer between two plates}

\subsection{Problem statement}

Consider a thin circular cylindrical layer $\Omega$ of thickness $2h$ and radius $a$. The layer thinness is represented by the inequality 
\begin{equation}\label{xi}
\xi:=\frac{h}{a} \ll 1 \:.
\end{equation}
The layer is made of a linear elastic material characterized by Lam\'e's constants $\mu$ and $\lambda$. Cylindrical coordinates for $\Omega$ are chosen so that
\begin{equation}\label{geometric definition}
r<a\:, \quad -h<z<h\:, \quad 0\leq \theta < 2\pi\:.
\end{equation}
We denote the top and bottom surface of $\Omega$ by $\partial\Omega^\pm$, and the cylindrical surface by $\partial\Omega^0$. 
The layer is perfectly bonded to two rigid plates, one is above $\partial\Omega^{+}$ and the other is below $\partial\Omega^{-}$, whereas $\partial\Omega^0$ is neither loaded nor constrained (Fig.~\ref{Fig1}a).

\bigskip
\noindent
We are interested in analyzing the linear elastic response of $\Omega$ when the plates are pulled apart, so that the top (bottom) plate is displaced by $U$ $(-U)$ along the $z$-axis. We calculate this response by solving an axisymmetric boundary-value problem of classical linear elasticity formulated for Navier's equations
\begin{eqnarray}
\left ( {\lambda}+ {\mu} \right)        \frac{\partial}{\partial r}\left( \frac{\partial u_r }{\partial r}+\frac{ u_r }{  r}+\frac{\partial u_z }{\partial z}\right)+\mu \left( \frac{\partial^2 u_r }{\partial r^2}+\frac{ 1 }{  r}\frac{\partial u_r }{\partial r}
-\frac{u_r}{r^2}+\frac{\partial^2 u_r }{\partial z^2}\right)=0\:, \label{equiulibrium r} \\
\left ( {\lambda}+ {\mu} \right)        \frac{\partial}{\partial z}\left( \frac{\partial u_r }{\partial r}+\frac{ u_r }{  r}+\frac{\partial u_z }{\partial z}\right)+\mu\left( \frac{\partial^2 u_z }{\partial r^2}+\frac{ 1 }{  r}\frac{\partial u_z }{\partial r}
+\frac{\partial^2 u_z }{\partial z^2}\right)=0 \:, \label{equiulibrium z}
\end{eqnarray}
with the boundary conditions on $\partial\Omega^{\pm}$ 
\begin{equation}\label{BC bottom and top}
u_z = \pm U\:, \quad u_r = 0\: ,
\end{equation}
and on $\partial\Omega^0$
\begin{equation}
\sigma_{rr} = \lambda \left( \frac{\partial u_r }{\partial r}+\frac{ u_r }{  r}+\frac{\partial u_z }{\partial z}\right) +2\mu \frac{\partial u_r }{\partial r} =0 \:, \label{normal cylindrical surface}
\end{equation}
\begin{equation}
\sigma_{rz}  =\mu \left( \frac{\partial u_r }{\partial z}+\frac{\partial u_z }{\partial r}\right) =0\label{shear cylindrical surface}\:.
\end{equation}

\subsection{Asymptotic solution}
We seek the solution using Love-Galerkin's bi-harmonic potential $\Phi$ \cite{Love,Galerkin}, so that the displacements are expressed as 
\begin{equation}\label{ur}
u_r = - \frac{1}{2(1-\nu)}  \frac{\partial^2 \Phi}{\partial r \partial z}
\end{equation}
and
\begin{equation}\label{uz}
u_z = \frac{\partial^2 \Phi}{\partial r^2}+ \frac{1}{r}\frac{\partial  \Phi}{\partial r }+ \frac{1-2\nu}{2(1-\nu)}\frac{\partial^2 \Phi}{\partial z^2}\:,
\end{equation} \\
and the potential must satisfy the bi-harmonic equation,
\begin{equation}\label{bi-harmonic}
\left(\frac{\partial^2 }{\partial r^2}+ \frac{1}{r}\frac{\partial }{\partial r }+ \frac{\partial^2 }{\partial z^2}\right)^2 \Phi =0 \:.
\end{equation}

\noindent
Since we are interested in thin layers, it is meaningful to introduce scaled dimensionless coordinates
\begin{equation}\label{scaled coordinates}
R: = \frac{r}{a} \quad {\rm and} \quad Z:= \frac{z}{h}= \frac{z}{\xi a }\:.
\end{equation}
Also in lieu of Poisson's ratio we use the dimensionless material parameter 
\begin{equation}\label{chi}
\chi : =\sqrt{ \frac{3(1-2\nu)}{2(1-\nu)}      }\:.
\end{equation}
The rationale behind this choice will become clear later, but for now it is sufficient to recognize that $\chi \rightarrow 0$ as $\nu \rightarrow 1/2$, and $0 \leq \chi \leq 3/2 $ as $1/2 \geq \nu \geq -1$. Thus nearly incompressible materials are characterized by $\chi \ll 1$. 

\bigskip
\noindent
Let us evaluate the bi-harmonic operator in the scaled coordinates:
\begin{equation}\label{bi-harmonic scaled}
\begin{split}
\left(\frac{\partial^2 }{\partial r^2}+ \frac{1}{r}\frac{\partial }{\partial r }+ \frac{\partial^2 }{\partial z^2}\right)^2 =
\frac{1}{a^2}\left(\frac{\partial^2 }{\partial R^2}+ \frac{1}{R}\frac{\partial }{\partial R }+ \frac{1}{\xi^2}\frac{\partial^2 }{\partial Z^2}\right)^2  \:.
\end{split}
\end{equation}
This expression implies that the last term is dominant and to the leading order the bi-harmonic equation is reduced to
$$
\frac{\partial^4 \Phi }{\partial Z^4 }=0\:,
$$
and therefore 
\begin{equation}\label{Love-Galerkin Potential}
\Phi =\xi^2 a^2 U \left\{  A_0 (R) +   A_1 (R) Z + A_2(R) Z^2 + A (R) Z^3 \right\} +{\cal O}(\xi^4) \:.
\end{equation}
Here $A$'s are unknown functions of $R$ to be determined by satisfying the boundary conditions. The pre-multiplier is dictated by elementary dimensional considerations, but of course its use is not essential. The term ${\cal O}(\xi^4)$  implies that we are content with the leading order term.

\bigskip
\noindent
Now we evaluate (\ref{ur}) and (\ref{uz}) in terms of (\ref{Love-Galerkin Potential}),  and substitute the displacements $u_r$ and $u_z$ in (\ref{BC bottom and top}). As a result we obtain
\begin{equation}\label{eq for the coefficients}
A_0(R)=0\:, \quad A_1(R)=-3A(R)\:,\quad {\rm and}\quad A_2 (R) = 0	\:,
\end{equation}
and the governing differential equation for $A (R)$:
\begin{equation}\label{equation for A}
A''(R)+ \frac{  
	1 }{ R}A'(R)
-
\left(\frac{\chi}{\xi}\right)^2 A (R) = -\frac{ 1}{2\xi^2 }\:.
\end{equation}
This equation can be solved exactly using Bessel functions:
\begin{equation}\label{equation for C3 intermediate} 
A (R) = 
C J_0\left(\frac{i  
	\chi }{\xi
}R\right)+C_1
Y_0\left(-\frac{i \chi
}{\xi }R\right)+\frac{1}{2 \chi ^2}\:.
\end{equation} \\
The first two terms on the right-hand side of this expression represent the homogeneous solution, and they involve two yet undetermined constants. Further, the simple dependence of the Bessel functions on the parameters $\xi$ and $\chi$ in (\ref{equation for C3 intermediate}) justifies the choice of $\chi$ defined in (\ref{chi}). The constant $C_1$ must be set equal to zero because $Y_0 (-ix)\rightarrow\infty$ as $x\rightarrow 0^+$. The constant $C$ is determined by satisfying the boundary conditions in (\ref{normal cylindrical surface}) and (\ref{shear cylindrical surface}). In the adopted setting, these boundary conditions cannot be satisfied exactly, but only in a Saint-Venant's sense,  for the normal and shear resultants:
\begin{equation}\label{normal BC resultant}
\int_{-1}^{1}{ \sigma_{rr}}_{|R=1} \: {\rm d}Z = 0
\end{equation}
and
\begin{equation}\label{shear BC resultant}
\int_{-1}^{1}{\sigma_{rz}}_{|R=1} \: {\rm d}Z = 0\:.
\end{equation} \\
The first of these conditions yields
\begin{equation}\label{equation for c1 intermediate} 
C = -\frac{3  \left(3-2 \chi
	^2\right)}{2 \chi ^2
	\left(3-\chi ^2\right)
	\left[3 I_0\left(\frac{\chi
	}{\xi }\right)-2 \xi  \chi 
	I_1\left(\frac{\chi }{\xi
	}\right)\right]} \:,
\end{equation}
whereas the second condition is trivially satisfied. 

\bigskip
\noindent
At this stage, we can use back substitution to obtain the expressions for Love-Galerkin's potential 
\begin{equation}\label{potential solution} 
\Phi = \frac{a^2 \xi^2
	U}{2 \chi ^2}
\left[  1 -\frac{ 3 \left(3-2 \chi
	^2\right)}{
	\left(3-\chi ^2\right)
	\left[ 3 I_0\left(\frac{\chi
	}{\xi }\right)-2 \xi  \chi 
	I_1\left(\frac{\chi }{\xi
	}\right)\right]} J_0\left(\frac{i  
	\chi }{\xi
}R\right)\right](-3Z+Z^3 )\:,
\end{equation}
and the displacements
\begin{equation}\label{ur solution} 
u_r = -\frac{3 \left(3- 2 \chi
	^2 \right)
	I_1\left(\frac{ \chi }{\xi
	}R\right)\left(1-Z^2 \right) U}{2 \chi \left[3 
	I_0\left(\frac{\chi }{\xi
	}\right)-2 \xi  \chi  
	I_1\left(\frac{\chi }{\xi
	}\right) \right] } \:,
\end{equation}
\begin{equation}\label{uz solution} 
u_z = \frac{  \left[ 6 \left(3-\chi
	^2\right)
	I_0\left(\frac{\chi }{\xi
	}\right)-4 \xi  \chi 
	\left(3-\chi ^2 \right)
	I_1\left(\frac{\chi }{\xi
	}\right)+3 \left(3- 2 \chi
	^2 \right)
	I_0\left(\frac{  \chi }{\xi
	}R\right)\left(1 - Z^2 \right)\right]ZU}{2
	\left(3-\chi ^2\right)
	\left[3
	I_0\left(\frac{\chi }{\xi
	}\right) - 2 \xi  \chi 
	I_1\left(\frac{\chi }{\xi
	}\right) \right]}\:.
\end{equation}

\bigskip
\noindent
The force acting on the upper (or lower) plate can be calculated as
\begin{equation}\label{F} 
\begin{split}
F &=    2\pi a^2 \int_0^1 { \sigma_{zz}}_{|Z= 1} R \: {\rm d}R= \frac{3 \pi   \mu aU}{8\xi^3}\times \\
&\frac{   8 \xi^2}{ \chi^3\left(3-\chi
	^2\right)} \left[
\frac{3 \chi  \left(3-\chi ^2\right)
	I_0\left(\frac{\chi }{\xi
	}\right)-2 \xi  \left(-\chi
	^4-3 \chi ^2+9\right)
	I_1\left(\frac{\chi }{\xi
	}\right) }{ 3
	I_0\left(\frac{\chi }{\xi
	}\right)-2 \xi  \chi 
	I_1\left(\frac{\chi }{\xi
	}\right)}
\right]\:.
\end{split}
\end{equation}
The rationale behind the particular factorization adopted in this equation will become clear in Section 2.3. Let us emphasize that (\ref{potential solution}-\ref{F}) are valid for any $\chi$.

\subsection{Particular cases}

In this section, we examine two limits for the solution given by (\ref{ur solution}-\ref{F}). First, in the limit as $\chi\rightarrow 0$ while $\xi$ is fixed, (\ref{ur solution}-\ref{F}) recover the classical solution for incompressible layers:
\begin{equation}\label{potential incompressible limit}
\Phi = -\frac{1}{8} a^2
\left(1-R^2\right)Z
\left(3-Z^2\right)U\:,
\end{equation}
\begin{equation}\label{ur incompressible limit} 
u_r = - \frac{3 R 
	\left(1-Z^2 \right)}{4 \xi } U \:,
\end{equation}
\begin{equation}\label{uz incompressible limit} 
u_z = \frac{1}{2}  Z
\left(3-Z^2\right) U\:,
\end{equation}
and
\begin{equation}\label{force incompressible limit} 
F =  \frac{3 \pi  \mu a U}{8
	\xi ^3}\:.
\end{equation} \\
Note that (\ref{ur incompressible limit}) and (\ref{uz incompressible limit}) coincide with (\ref{fluids radial velocity}) and (\ref{fluids axial velocity}), respectively, and (\ref{F}) is expressed in a way that results in the function displayed on the second line to become unity in the limit as $\chi\rightarrow 0$ while $\xi$ is fixed. If needed, one can calculate the deviatoric components from (\ref{ur incompressible limit}) and (\ref{uz incompressible limit}), and then determine the pressure from the equilibrium equations. 

\begin{figure}
	\begin{center} 
		\includegraphics[width=\textwidth]{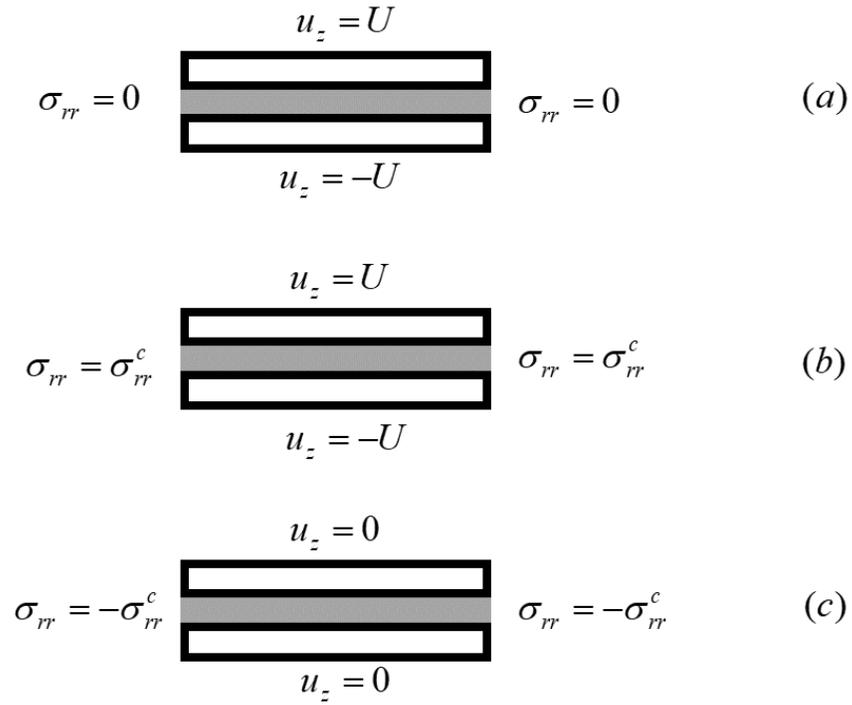}
		\caption{The problem for a compressible layer (a) as a superposition of uniaxial straining (b) and a boundary-layer  problem (c). The stress $\sigma_{rr}^c$ is $\sigma_{rr}$ from (\ref{BC cylindrical surface 2}).}		\label{Fig2} 
	\end{center}
\end{figure}

\bigskip
\noindent
The second particular case represents compressible layers characterized by $\chi = {\cal O}(1)$. To evaluate this case we introduce the variable 
\begin{equation}\label{zeta}
\zeta : =\frac{\xi}{\chi}\:,
\end{equation}
so that compressible layers can be associated with the limit $\zeta \rightarrow 0$ as $\xi$ is fixed. In evaluating this limit, we use
$$
\lim_{\zeta\rightarrow 0}  \frac{ \zeta I_1(  \zeta^{-1} )}{ I_0(   \zeta^{-1} )} =0\:,
$$
to obtain
\begin{equation}\label{potential compressible limit} 
\Phi=\frac{1}{2}a^2 \left( \frac{\xi}{\chi} \right)^2 Z^3 U \:,
\end{equation}
\begin{equation}\label{ur compressible limit} 
u_r =0 \:,
\end{equation}
\begin{equation}\label{uz compressible limit} 
u_z = Z U\:,
\end{equation}
and
\begin{equation}\label{force compressible limit} 
F =  \frac{3 \pi  \mu a U}{
	\xi \chi^2}\:.
\end{equation} \\
These displacement fields and force represent the case of uniaxial straining along the $z$-axis, so that the only non-zero strain component is
$$ \epsilon_{zz}=\frac{U}{h}=\frac{U}{a\xi} \:.$$
The corresponding stress state is
\begin{equation}\label{compressible stress state}
\sigma_{rr}=\sigma_{\theta\theta}=\frac{ \mu  U  \left(3-2
	\chi ^2\right)}{a \xi  \chi
	^2}\:,\quad\quad\sigma_{zz}=\frac{3 \mu  U}{a \xi  \chi
	^2}\:,\quad\quad \sigma_{rz}=0\:.
\end{equation}
This stress state has two remarkable properties. First, the shear stress on the plate-layer interfaces is equal to zero, so that this stress state does not differentiate between no-slip and shear-traction-free conditions on the plate-layer interfaces. Second, the stress state results in non-zero normal traction component on the cylindrical surface $\partial\Omega_0$, and no averaging procedure can change that simply because the stress $\sigma_{rr}$ is constant. Thus the adopted basic asymptotic ansatz does not work in the limit as $\zeta \rightarrow 0$ while $\xi$ is fixed. This issue can be addressed by reformulating the original boundary-value as a superposition of two problems (Fig.~\ref{Fig2}). The first problem is identical to the original problem, except for the new boundary condition on $\partial\Omega_0$ replacing (\ref{normal cylindrical surface}) and (\ref{shear cylindrical surface}):
\begin{equation}\label{BC cylindrical surface 1}
\sigma_{rr} = \frac{ \mu  U  \left(3-2
	\chi ^2\right)}{a \xi  \chi
	^2}\:, \quad \sigma_{rz}  = 0 \quad {\rm on}\:\partial\Omega^0\:.
\end{equation} \\
The second problem is characterized by the boundary conditions
\begin{equation}\label{BC bottom and top 2}
u_z = 0\:, \quad u_r = 0 \quad {\rm on}\:\: \partial\Omega^{\pm}\:,
\end{equation}
and
\begin{equation}\label{BC cylindrical surface 2}
\sigma_{rr} = -\frac{ \mu  U  \left(3-2
	\chi ^2\right)}{a \xi  \chi
	^2}\:, \quad \sigma_{rz}  = 0 \quad {\rm on}\:\: \partial\Omega^0\:.
\end{equation} \smallskip

\noindent Then the solution for the first problem is given by (\ref{ur compressible limit}-\ref{force compressible limit}). The second-problem is of the boundary-layer type. It is characterized by an exponential decay away from the cylindrical surface, but the decay length may significantly depend on the problem parameters \cite{Maz'ya, Movchans1}. This issue will be examined in the next section, where we compare asymptotic and finite element solutions. 

\subsection{Comparisons with finite element results}

The problem of interest is straightforward to analyze using a finite element method, and we did it using the commercial program ABAQUS. We  employed uniform meshes formed by square eight-node hybrid elements; CAX8H in the ABAQUS language. This choice of element type was beneficial for analyzing problems in a wide range of $\chi$. The mesh size was chosen based on convergence studies, so that the force $F$ was computed accurately up to the first six significant digits.  Comparisons between the asymptotic and finite element solutions for the force are presented in Table~\ref{error for our solution for the plates}, where the finite element solutions are treated as the exact ones.   Clearly, the asymptotic solution (\ref{F}) for the force becomes progressively more accurate as $\xi$ decreases in the entire range of $\chi$. The number of correct significant digits increases by one with each order of magnitude decrease in $\xi$. In particular, for $\xi=10^{-3}$, the asymptotic solution accurately captures the first three significant digits across all $\chi$. 
\begin{table}[htb]
	\begin{center}
		\begin{tabular}{| c | c | c | c |}
			\hline
			$\quad$ & $ \: \xi = 10^{-3} $ & $ \: \xi = 10^{-2} $ & $ \: \xi = 10^{-1} $ \Tstrut\Bstrut\\ \hline
			$ \: \chi = 10^{-3} $ & $ \: 3 \times 10^{-4} $ & $ \: 3 \times 10^{-3} $ & $ \: 5 \times 10^{-2}$  \Tstrut\Bstrut\\ \hline
			$\: \chi = 10^{-2} $ & $ \: 2 \times 10^{-4} $ & $ \: 3 \times 10^{-3} $ & $ \: 5 \times 10^{-2}$  \Tstrut\Bstrut\\ \hline
			$\: \chi = 10^{-1} $ & $ \: 2 \times 10^{-4} $ & $ \: 7 \times 10^{-4} $ & $ \: 5 \times 10^{-2}$  \Tstrut\Bstrut\\ \hline
			$ \: \chi = 1 $       & $ \: 8 \times 10^{-4} $ & $ \: 8 \times 10^{-3} $ & $ \: 8 \times 10^{-2} $  \Tstrut\Bstrut\\ \hline
		\end{tabular}
	\end{center}
	\caption{Errors of the asymptotic solution for the force (\ref{F}) in comparison to finite element solutions accurate up to the first six significant digits.}\label{error for our solution for the plates}
\end{table}

\begin{figure}[h]
	\begin{center} 
		\includegraphics[width=\textwidth]{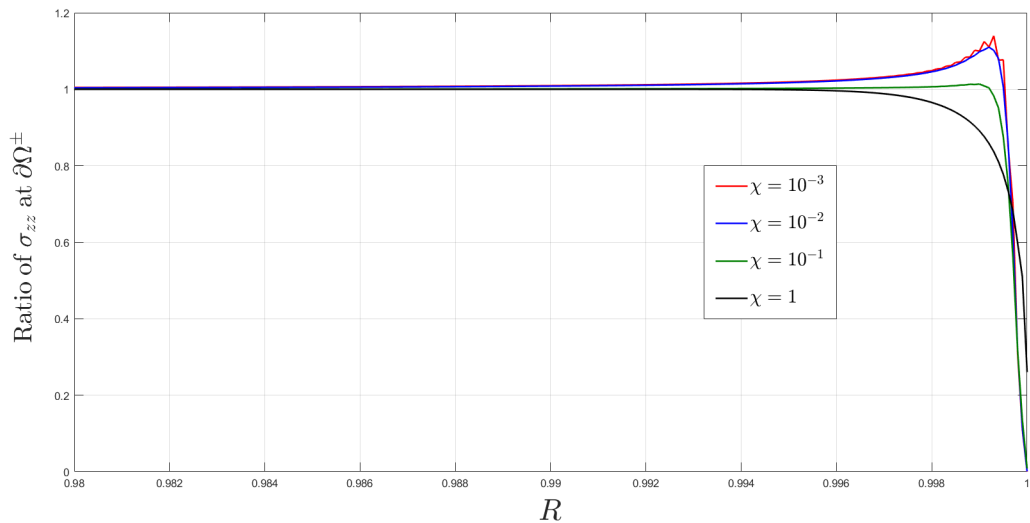}
		\caption{The ratio of the asymptotic and finite element solutions for  $\sigma_{zz}$ along the interface at $Z =  1$  for $\xi = 10^{-3}$.}		\label{Fig3} 
	\end{center}
\end{figure}

\begin{figure}[h]
	\begin{center} 
		\includegraphics[width=\textwidth]{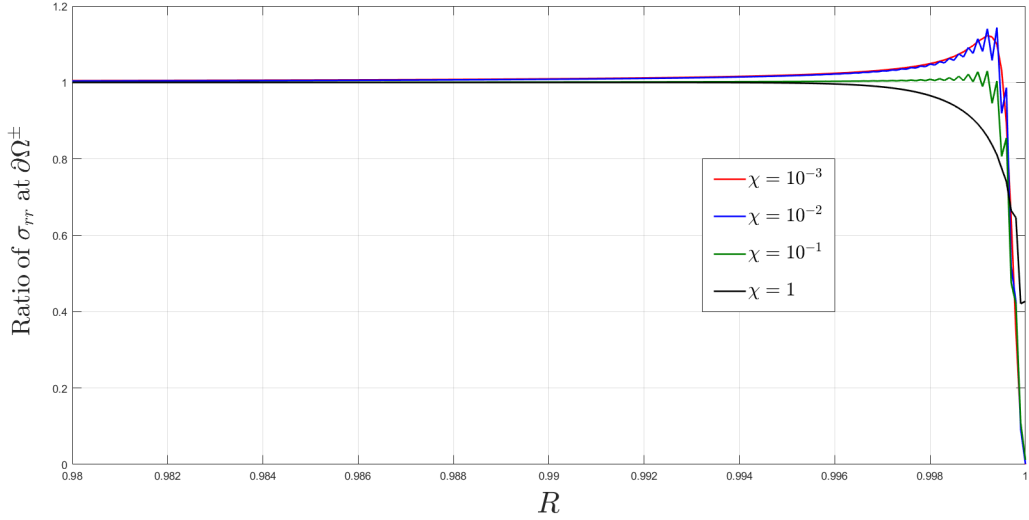}
		\caption{The ratio of the asymptotic and finite element solutions for  $\sigma_{rr}$ along the interface at $Z =  1$  for $\xi = 10^{-3}$.}		\label{Fig4} 
	\end{center}
\end{figure}

\noindent
The fact that the asymptotic solution for the force agrees well with the finite element solutions implies that the boundary-layer correction is insignificant, at least as far as the force calculations are concerned. To get an additional confirmation of this statement we compared the asymptotic versus finite element solutions for the stress $\sigma_{zz}$ along the interface  $Z=  1$.  To this end, we plotted the ratio of the asymptotic versus finite element solution as a function of $R$ for 
$\xi=10^{-3}$ and $\chi = 10^{-3},10^{-2},10^{-1},1$ (Fig.~\ref{Fig3}). Note that we used $0.98\leq R \leq 1$ for plotting because the ratio is very close to unity for $R<0.98$. It is clear that the asymptotic solution is incapable of predicting singular stresses near the corner, as it does not include the boundary-layer asymptotic corrections. Those corrections decay exponentially, and the rate of decay mildly depends on $\chi$. Apparently the region of dominance of those singular stresses is sufficiently small, so that the singularity does not significantly affect the force. 
We further examined the stress field by plotting the stress $\sigma_{rr}$ in the same manner as we did for $\sigma_{zz}$ (Fig.~\ref{Fig4}). In this case, the objective was to examine Saint-Venant's effect associated with the weak imposition of the boundary conditions at the surface $\partial\Omega_0$. Figure~\ref{Fig4} suggests that the strong boundary-layer effect decays exponentially fast over distances comparable to $\xi=10^{-3}$, which is consistent with Saint-Venant's principle.
The observed exponential decays are consistent with the analysis of \cite{Maz'ya}, for the case of Dirichlet boundary conditions on the flat surfaces. While the decay rate mildly depends on  $\chi$, which makes the analysis interesting, there are no surprises here -- the boundary layer is localized and has no significant effect on the force.

\subsection{Comparisons with other approximate solutions}
Due to significance of stretching/squeezing of thin layers between stiff plates for mechanical characterization of rubber and other polymers, there is a large number of approximate solutions available in the literature.
As far as testing is concerned, the most accessible specimen property is the so-called apparent modulus defined as the ratio of the average axial stress and strain
\begin{equation}\label{E_a}
E_A := \frac{F}{\pi a^2 } \div \frac{U}{h} \: . 
\end{equation}
It is natural to normalize the apparent modulus with the actual Young's modulus of the layer, and define 
\begin{equation}\label{key}
\hat{E}: =\frac{ E_A }{E}\:.
\end{equation}
The normalized apparent modulus according to our analysis follows directly from (\ref{F}):
\begin{equation}\label{YA}
\hat{E} = \frac{9 \chi  \left(3-\chi ^2\right) I_0\left(\frac{\chi }{\xi }\right)-6 \xi  \left(9-3 \chi^2-\chi^4\right) I_1\left(\frac{\chi }{\xi }\right)}{\chi ^3 \left(9-4 \chi ^2\right) \left[3
	I_0\left(\frac{\chi }{\xi }\right)-2 \xi  \chi  I_1\left(\frac{\chi }{\xi }\right)\right]} \: .
\end{equation} 
The two extremes of this expression are 
\begin{equation}\label{AY incomp}
\hat{E}^i = \frac{1}{8 \xi ^2} \quad {\rm for } \: \chi=0\: ,
\end{equation}
and 
\begin{equation}\label{AY comp}
\hat{E}^c = \frac{3 \left(3-\chi ^2\right)}{\chi ^2 \left(9-4 \chi ^2\right)} \quad {\rm for } \: \chi={\cal O} (1) \: . 
\end{equation} \\   
These expressions can be derived either directly from (\ref{YA}) or from (\ref{force incompressible limit}) and  (\ref{force compressible limit}).  

\bigskip
\noindent
Since historically the configuration was used for testing rubber, it is appropriate to compare various approximate solutions by choosing $\nu =0.499905$ \cite{Anderson2004++}, which corresponds to $\chi=0.0238724$. 
In Table~\ref{apparent moduli for the plates} we present relative errors  for various approximate solutions for $\hat{E}$, when compared to finite element solutions accurate up to the first six significant digits. The approximate solutions are arranged in the chronological order and the errors are computed for $\xi=10^{-1},10^{-2},10^{-3}$. 
\begin{table}[h]
	\begin{center}
		\centerline{\begin{tabular}{| c | c | c | c |}
				\hline
				$\quad$ & $\: ~ \xi = 10^{-3} \: $ & $\: ~ \xi = 10^{-2} \: $ & $ \: ~ \xi = 10^{-1} $ \Tstrut\Bstrut\\ \hline
				$~$Ref. \cite{GL1959++}   & \cellcolor{gray!25} $ ~ 7 \times 10^{-2} \: $ & \cellcolor{gray!25} $ ~ 1 \times 10^{-1}\:$ & \cellcolor{gray!25}$ ~ 3 \times 10^{-2}  $\Tstrut\Bstrut\\ \hline
				$~$Ref. \cite{LS1963++}     & $ \: \: \: 7 \times 10^{-5} \: $ & $ \: \: \: 3 \times 10^{-3} \: $ & $~ 3 \times 10^{-2} $\Tstrut\Bstrut\\ \hline
				$~$Ref. \cite{L1979++}   &  $ \: \: \: 4 \times 10^{-3} \: $ & $ \: \: 5 \times 10^{-3} $ & $~ 3 \times 10^{-2} $\Tstrut\Bstrut\\ \hline
				$~$Ref. \cite{CK1990++} &$ \: \: 2 \times 10^{-4} $ & $ \: \: 2 \times 10^{-3} $ & $~ 5 \times 10^{-2} $\Tstrut\Bstrut\\ \hline
				$~$Ref. \cite{G1994++} & $\: \: 4 \times 10^{-4} $ & $ \: \: 3 \times 10^{-3} $ & $~ 3 \times 10^{-2} $\Tstrut\Bstrut\\ \hline
				$~$Ref.  \cite{TL1998++}         & $ \: \: 7 \times 10^{-5} $ & $ \: \: 3 \times 10^{-3} $  & $~ 3 \times 10^{-2} $\Tstrut\Bstrut\\ \hline
				$~$Ref.  \cite{H2002++}  & \cellcolor{gray!25} $ ~ 7 \times 10^{-2} \: $  & \cellcolor{gray!25} $~ 1 \times 10^{-1}\:$ & \cellcolor{gray!25}$ ~ 4 \times 10^{-2}$\Tstrut\Bstrut\\ \hline
				$~$Ref.  \cite{QL2015++}       & $ \: \: 6 \times 10^{-5} $ & $ \: \: 3 \times 10^{-3} $ & $ ~ 3 \times 10^{-2} $\Tstrut\Bstrut\\ \hline
					$~$Ref.  \cite{RAS2}    &$ \: \: 2 \times 10^{-4} $ & $ \: \: 2 \times 10^{-3} $ & $~ 5 \times 10^{-2} $\Tstrut\Bstrut\\ \hline
				$\:~$Eq.  (\ref{YA}) & $ \: \: 2 \times 10^{-5} $ & $ \: \: 2 \times 10^{-3} $ & $~ 5 \times 10^{-2} $\Tstrut\Bstrut\\ \hline 
		\end{tabular}}
	\end{center}
	\caption{Relative errors in the normalized apparent modulus of rubber when compared against finite element results for $\nu=0.499905$ \cite{Anderson2004++}. Rejected solutions are shown in gray.}
	\label{apparent moduli for the plates}
\end{table}
It is clear that, beside our solution, seven solutions stand out as very accurate: Lindsey's {\it et al.} \cite{LS1963++}, Lindley's \cite{L1979++}, Chalhoub and Kelly's \cite{CK1990++}, Gent's \cite{G1994++}, 
		Tsai and Lee's \cite{TL1998++}, Qiao and Lu's \cite{QL2015++}, and Schapery's \cite{RAS2}. All seven solutions relied on satisfying the boundary conditions on $\partial\Omega^0$ weakly, and therefore, not surprisingly,  all of them improve as
		 $\xi \rightarrow 0$. Further, all chosen solutions recover exactly the limit as $\chi \rightarrow 0$. In \cite{LS1963++,L1979++,CK1990++,G1994++, TL1998++,  RAS2} it was assumed that $u_r$ has a parabolic profile. This is an excellent assumption, confirmed by our asymptotic and finite element analyses. Qiao and Lu assumed that each displacement field  can be expressed as a separable product  of a function of $R$ and a function of $Z$.  This assumption holds for $u_r$, but, according to our analysis, $u_z$ requires two  products rather than one. The expressions derived by Lindsey {\it et al.} \cite{LS1963++}  and	Tsai and Lee  \cite{TL1998++} are identical, although the approaches are different. The same statement applies to the solutions of Chalhoub and Kelly \cite{CK1990++} and Schapery \cite{RAS2}.  We rejected two solutions because they are characterized by significant errors, and the errors  do not decay as $\xi$ decreases. Let us mention that concerns about \cite{H2002++} have been expressed in \cite{pence}.

\bigskip
\noindent
In the next verification round, we considered the eight approximations which passed the first round, and compared them with accurate finite element solutions for $\nu=0.3$. Results of those comparisons are summarized in Table~\ref{apparent moduli for plates round 2}.  Among those approximations three, shown in gray, were rejected. This was not surprising as the three rejects were constructed for nearly incompressible materials. 
\begin{table}[h]
	\begin{center}
		\centerline{\begin{tabular}{| c | c | c | c |}
				\hline
				$\quad$ & $\: ~ \xi = 10^{-3} \: $ & $\: ~ \xi = 10^{-2} \: $ & $ \: ~ \xi = 10^{-1} $ \Tstrut\Bstrut\\ \hline
				$~$Ref. \cite{LS1963++}     & $ \: \: \: 2 \times 10^{-5} \: $ & $ \: \: \: 1 \times 10^{-4} \: $ & $~ 2 \times 10^{-3} $\Tstrut\Bstrut\\ \hline
				$~$Ref. \cite{L1979++}   &  $ \: \: \: 3 \times 10^{-5} \: $ & $ \: \: 3 \times 10^{-4} $ & $~ 2 \times 10^{-3} $\Tstrut\Bstrut\\ \hline
				$~$Ref. \cite{CK1990++} &\cellcolor{gray!25}$ \: \: 4 \times 10^{-1} $ & \cellcolor{gray!25}$ \: \: 4 \times 10^{-1} $ & \cellcolor{gray!25}$~ 5 \times 10^{-1} $\Tstrut\Bstrut\\ \hline
				$~$Ref. \cite{G1994++} & \cellcolor{gray!25}$\: \: 4 \times 10^{-1} $ & \cellcolor{gray!25}$ \: \: 4 \times 10^{-1} $ & \cellcolor{gray!25}$~ 3 \times 10^{-1} $\Tstrut\Bstrut\\ \hline
				$~$Ref.  \cite{TL1998++}         & $ \: \: 2 \times 10^{-5} $ & $ \: \: 1 \times 10^{-4} $  & $~ 2 \times 10^{-3} $\Tstrut\Bstrut\\ \hline
				$~$Ref.  \cite{QL2015++}       & $ \: \: 3 \times 10^{-6} $ & $ \: \: 5 \times 10^{-5} $ & $ ~ 3 \times 10^{-4} $\Tstrut\Bstrut\\ \hline
				$~$Ref. \cite{RAS2} &\cellcolor{gray!25}$ \: \: 4 \times 10^{-1} $ & \cellcolor{gray!25}$ \: \: 4 \times 10^{-1} $ & \cellcolor{gray!25}$~ 5 \times 10^{-1} $\Tstrut\Bstrut\\ \hline
				$\:~$Eq.  (\ref{YA}) & $ \: \: 5 \times 10^{-4} $ & $ \: \: 5 \times 10^{-3} $ & $~ 5 \times 10^{-2} $\Tstrut\Bstrut\\ \hline
		\end{tabular}}
	\end{center}
	\caption{Relative errors in the normalized apparent modulus when compared against finite element results for $\nu=0.3$. Rejected solutions are  shown in gray.}
	\label{apparent moduli for plates round 2}
\end{table}

\noindent The remaining five were compared with each other for $0<\nu<0.49999$ and they were in excellent agreement. 	 

\bigskip
\noindent
The approximation obtained by  Lindsey {\it et al.}  \cite{LS1963++} and later by Tsai and Lee  \cite{TL1998++} has the form
\begin{equation}\label{YAS}
\hat{E}_L = \frac{\left(3-\chi ^2\right) \left[9 \chi  I_0\left(\frac{\chi }{\xi }\right)-2 \xi  \left(4 \chi ^4-9
	\chi ^2+9\right) I_1\left(\frac{\chi }{\xi }\right)\right]}{\chi ^3 \left(9-4 \chi ^2\right) \left[3
	I_0\left(\frac{\chi }{\xi }\right)-2 \xi  \chi  I_1\left(\frac{\chi }{\xi }\right)\right]} \: ,
\end{equation} 
which is very close to (\ref{YA}). Indeed the denominators of the two expressions coincide, and the difference between the two approximations is estimated as
\begin{equation}\label{YASError}
\frac{\hat{E}_L - \hat{E} }{ \hat{E} } \approx -\frac{2 \chi  \left(4 \chi
	^4-24 \chi ^2+27\right)}{9
	\left(3-\chi ^2\right)} \xi \: .
\end{equation} 
The function of $\chi$ in this estimate is equal zero at  $\chi = 0 $, as expected, and it varies between minus one and plus three. Thus we conclude that the difference between $ \hat{E}_L$ and  $\hat{E}  $ is small in the entire range of $\chi$. Analytical comparisons of our approximation  with those of Lindley's \cite{L1979++} and Qiao and Lu \cite{QL2015++} are somewhat cumbersome because the former has a conditional structure, and the latter is predicated on solving a transcendental equation. 

\bigskip
\noindent
    It is remarkable that the approximate solution in \cite{LS1963++}, presented more than fifty years ago, is accurate in the entire range of $\chi$.  That solution was developed by assuming the displacements in the form
	\begin{equation}\label{Sch1}
	u_r (R,Z) = f(R)(1-Z^2)
	\quad{\rm and}\quad
	u_z (R,Z) = UZ\:,
	\end{equation}
with the function $f(R)$  determined upon averaging the governing differential equations and boundary conditions on the cylindrical surface through the thickness. It turns out that $f(R)$ derived in this manner 
gives rise to $u_r$ which coincides with our solution (\ref{ur solution}). Of course $u_z$ in (\ref{Sch1}) is different from our solution (\ref{uz solution}) simply because our solution is a cubic polynomial in $Z$.
Recently, Schapery \cite{RAS2}, who is a co-author of \cite{LS1963++}, proposed  the displacements in the form
	\begin{equation}\label{Sch2}
u_r (R,Z) = f(R)(1-Z^2)
\quad{\rm and}\quad
u_z (R,Z) = UZ +g(R) (Z-Z^3)\:.
\end{equation}
He did not derive the function $g(R)$, but it can be chosen to match our solution (\ref{uz solution}) exactly. Further, Schapery  suggests that (\ref{Sch2}) can yield accurate approximations without requiring either $\xi$ or $\chi$ to be small. We disagree with this position, as we believe that it is essential for $\xi \ll 1$, because, without this assumption, Saint-Venant's principle becomes meaningless for the problem at hand. 

\bigskip
\noindent
Schapery did not derive the function  $g(R)$ because he was interested in nearly incompressible layers. In this case, the dominant contribution to the force (or apparent modulus) is the pressure derived from the equilibrium equation along $r$, whose leading order asymptotic form is
$$
- \frac{\partial p }{\partial R} +\frac{\mu}{\xi}\frac{\partial u_r }{\partial Z}=0\:.
$$
Thus, for nearly incompressible layers, it is important to accurately calculate $u_r$ but not $u_z$. This also explains why (\ref{Sch1}) is accurate at least for nearly incompressible layers. Of course the fact that (\ref{Sch1}) works for all $\chi$ is a tribute to the clever averaging scheme that yielded $f(R)$.

\bigskip
\noindent
The ansatz in (\ref{Sch2}) provides a natural connection with \cite{LS1963++} and our work, but it is not the central theme of \cite{RAS2}. Rather that work focuses on an asymptotic method restricted to nearly incompressible layers, which allows one to solve a large class of practically important three-dimensional problems \cite{RAS2,RAS1}. For this reason, Schapery's approximate solution compared well with finite element results for nearly incompressible but not compressible cases. We will discuss Schapery's approach to three-dimensional problems for nearly incompressible layers later, in connection to an alternative asymptotic methodology, applied to analysis of thin layers between spheres.

\section{Thin layer between two spheres}

The problem for a thin layer between two equal rigid spheres of radius $a$ and minimum separation distance $2h$ (Fig.~\ref{Fig1}b) can be formulated and solved similarly to the problem for two plates. The mathematical differences turned out to be mostly technical rather than conceptual, as the equations become complicated to an extent that one should heavily rely on symbolic manipulators. To this end, let us mention that the ordinary differential equation which parallels (\ref{equation for A}) could only be solved with Maple but not  Mathematica. Due to similarities between the two problems and very long expressions for the case of two spheres, our presentation focuses on emphasizing the differences between the two problems rather than solution details.

\subsection{Problem statement}
In cylindrical coordinates shown in Figure~\ref{Fig1}b, the spheres are  prescribed by the equations
\begin{equation}\label{sphere physical coordinates upper}
r^2 + (z -  h - a)^2 = a^2 \quad(\text{upper sphere}) 
\end{equation}
and
\begin{equation}\label{sphere physical coordinates lower}
r^2 + (z +  h + a)^2 = a^2 \quad(\text{lower sphere}). 
\end{equation} \\
We assume that the layer is bounded by a circular cylindrical surface $\partial\Omega^0$ of radius $a$. Then the bounding surfaces 
 $\partial\Omega^+$ (top) and $\partial\Omega^-$ (bottom) are hemi-spheres prescribed by (\ref{sphere physical coordinates upper}) and (\ref{sphere physical coordinates lower}), respectively. 
 With these definitions, the boundary-value problem for a thin layer between two equal rigid spheres is prescribed by (\ref{equiulibrium r}-\ref{shear cylindrical surface}), with the provision that 
 (\ref{normal cylindrical surface}) and (\ref{shear cylindrical surface}) are imposed weakly, similar to (\ref{normal BC resultant}) and (\ref{shear BC resultant}).

\subsection{Asymptotic solution}
Following \cite{Jeffrey} we introduce the scaled coordinates, 
\begin{equation}\label{scaled coordinates spheres}
R: = \frac{r}{a\sqrt{\xi}} \quad {\rm and} \quad Z:= \frac{z}{h}= \frac{z}{\xi a }\:,
\end{equation}
different from those adopted in the previous section.
In these coordinates, to a leading order,  (\ref{sphere physical coordinates upper}) and (\ref{sphere physical coordinates lower}) can be combined in  the form
\begin{equation}\label{sphere scaled coordinates}
Z =\pm \left( 1  + \frac{1}{2}  R^2\right)\:. 
\end{equation}
The new scaled coordinates do not affect the leading order approximation for the bi-harmonic operator, but Love-Galerkin's potential has a slightly different form, 
\begin{equation}\label{Love-Galerkin Potential spheres}
\Phi =\xi a^2 U \left\{   A_0 (R) +  \left[ \int  A_1 (R) \: {\rm d}R \right] Z + A_2(R) Z^2 + A (R) Z^3 \right\} +{\cal O} \left(\xi^2\right)\:.
\end{equation}
This form differs from (\ref{Love-Galerkin Potential}) in two aspects. First, the pre-multiplier in (\ref{Love-Galerkin Potential spheres}) is consistent with the scaled coordinates in (\ref{scaled coordinates spheres}) rather than (\ref{scaled coordinates}). Second, $A_1 (R)$ is replaced with its anti-derivative to simplify the calculations.

\bigskip
\noindent 
The equations for the $A$-functions in (\ref{Love-Galerkin Potential spheres}) are
\begin{equation}\label{eq for the coefficients spheres}
A_0(R)=0\:, \quad A_1  (R)=-3
\left(1 + \frac{1}{2}R^2\right)^2 A'(R)\:,\quad {\rm and}\quad A_2 (R) = 0	\:,
\end{equation}
and 
\begin{equation}\label{equation for A spheres}
A''(R)+ \frac{7 R^2+2}{R^3+2 R} A'(R) -\frac{4\chi^2 }{\xi
	\left(R^2+2\right)^2}  A(R)=-\frac{4}{\left(R^2+2\right)^3}\:.
\end{equation} \\
This differential equation can be solved in terms of hypergeometric functions using Maple but not Mathematica:
\begin{eqnarray} \label{maple solution}
A(R)&=&  C_1 
\left(1 + \frac{1}{2}R^2  \right)^{-1-\beta} \, _2F_1\left(-1 -\beta
 ,2-\beta, 1-2 \beta
,1+ \frac{1}{2}R^2 \right) \nonumber \\ &+& C_2
\left(1 + \frac{1}{2}R^2  \right)^{-1+ \beta }
\, _2F_1\left(-1+ \beta  ,2+\beta
, 1+ 2\beta,
 1+ \frac{1}{2}R^2 \right)\nonumber \\
&+&\frac{1+ \beta ^2+R^2 }{4 \beta ^2
	\left(1- \beta ^2 \right)
	\left(1 + \frac{1}{2}R^2  \right)}\:, 
\end{eqnarray}
with
\begin{equation}\label{beta}
\beta = \sqrt{ 1 -  \frac{ \chi^2 }{ 2 \xi } }\:.
\end{equation}
The integration constants $C_1$ and $C_2$ can be determined using the conditions that the solution must be finite at $R=0$,  and the resultants on the cylindrical surface must be equal to zero:
$$
\int_{-1-1/2/\xi}^{1+1/2/\xi} {\sigma_{rr} }_{|R=1/\sqrt{\xi}} \: {\rm d}Z = \int_{-1-1/2/\xi}^{1+1/2/\xi} {\sigma_{rz} }_{|R=1/\sqrt{\xi}} \: {\rm d}Z = 0\:.
$$
As in the case of flat layers, the last condition is trivially satisfied and the other two result in 
\begin{eqnarray} \label{our solution}
A(R)&=&  C  \Gamma (-\beta -1)
\Gamma (2-\beta )
\left(1+ \frac{1}{2}R^2 \right)^{-1-\beta  }
\, _2{F}_1\left(-1 -\beta , 2-\beta , 1-2 \beta, 1+ \frac{1}{2}R^2\right)\nonumber \\
&-& C \frac{\Gamma (-1 + \beta)
\Gamma (2+\beta)}{\Gamma (1+ 2 \beta)}
\left(1+ \frac{1}{2}R^2 \right)^{-1+\beta }
\,  _2F_1\left(-1+\beta,2+\beta,1+ 2 \beta , 1+ \frac{1}{2}R^2\right) \nonumber 
\\
&+&\frac{1+ \beta^2+R^2 }
{4 \beta ^2
	\left(1- \beta ^2 \right)
	\left(1 + \frac{1}{2}R^2  \right)}\:,
\end{eqnarray}
with
\begin{eqnarray*}
C&=&\frac{3 
	(1-\beta )	\left(1+\beta^2\right)\Gamma (2-\beta ) \Gamma
	(-\beta )}{2^{ 1- 3\beta} \beta^2   \xi^{ \beta} }\times   \\
&&\left[ 6
\left(1+\beta \right)
\,
_2\tilde{F}_1\left(3-\beta
,-\beta ,2-2 \beta
,\frac{1}{2 \xi }\right) 
-
(3-\beta ) (4+ \beta )
\,
_2\tilde{F}_1\left(4-\beta
,-\beta ,2-2 \beta
,\frac{1}{2 \xi }\right) \right]\\
&+& \frac{3   
	\left(1+\beta \right)	\left(1+\beta^2\right)\Gamma
	(\beta ) \Gamma (2+\beta)}{2^{ 1-\beta} \beta^2 \xi^{3\beta} }\times\\
&&\left[       3  (2+ \beta  )
 \, _2\tilde{F}_1\left(\beta
,3+\beta ,2 +2 \beta
, \frac{1}{2 \xi
}\right)
+\beta (4-\beta ) 
\, _2\tilde{F}_1\left(1+ \beta
,3+\beta , 2+2 \beta
,\frac{1}{2 \xi
}\right)       \right]\:.
\end{eqnarray*}
In these equations, the tilde denotes the regularized hypergeometric function and $\Gamma$ is Euler's $\Gamma$-function. Let us mention that $\beta$ can take on both real and imaginary values, and the hypergeometric functions can take on complex values. Nevertheless, once the boundary conditions have been imposed, the resulting function $A(R)$ is real-valued.

\bigskip
\noindent
It is clear that the displacement and stress  fields are straightforward to derive  from (\ref{our solution}) by differentiation using a symbolic manipulator, but the resulting expressions are too long and hardly provide any insight. The force, 
\begin{equation} \label{F spheres}
\begin{split}
F = 2\pi a^2\xi \int_{0}^{1/\sqrt{\xi}} {\sigma_{zz}}_{|Z=1+\frac{1}{2}R^2} \: R \: {\rm d } R = 6\pi a \mu U \Psi ( \chi , \xi )\:,
\end{split}
\end{equation} 
can be calculated only numerically. 
Here $\Psi ( \chi , \xi )$ is introduced as a dimensionless force, whose behavior at the extremes of $\chi=0$ and $\chi={\cal O} (1)$ are known from \cite{Jeffrey} and \cite{PTK}, respectively:
\begin{equation} \label{incompressible F spheres}
\Psi^i ( 0 , \xi ) = \frac{1 }{4 \xi}
\end{equation} 
and
\begin{equation}\label{compressible F spheres}
\Psi^c ( \chi , \xi ) =  \frac{1}{\chi^2} \log{ \frac{1}{2\xi}} 
\:.
\end{equation}

\bigskip
\noindent
Rather than verifying $A(R)$ by establishing that the corresponding fields recover their counterparts provided in \cite{PTK} and  \cite{Jeffrey}, we chose to compare $\Psi ( \chi , \xi ) $ directly with  with finite element solutions accurate up to the first six significant digits for $\xi = 10^{-5},  10^{-4},  10^{-3},  10^{-2}$ and $\chi =  10^{-3},  10^{-2},  10^{-1}, 1$. This is done in Table~\ref{forces for the spheres}, where we also include the predictions based on  $\Psi^i (0 , \xi ) $, and $\Psi^c ( \chi , \xi ) $. Thus for each pair of  $\xi$ and $\chi$ we compare four numbers. The data clearly shows that $\Psi^i ( 0 , \xi )$ and $\Psi^c ( \chi , \xi )$ match the finite element solutions as expected: $\Psi^i ( 0 , \xi )$ works well for $\chi =  10^{-3}$ and $\Psi^c ( \chi , \xi )$ works well for $\chi =  1$. Both $\Psi^i (0, \xi ) $, and $\Psi^c ( \chi , \xi ) $ are outright inapplicable for $\chi =  10^{-1}.$ Note that $\Psi^i ( 0 , \xi )$ works well for $\chi =  10^{-2}$ for all $\xi$ except for $\xi =  10^{-5}$, which is somewhat perplexing because $\Psi^i ( 0 , \xi )$ is expected to be increasingly accurate as $\xi$ decreases. In contrast $\Psi ( \chi , \xi ) $ corresponding to (\ref{our solution}) matches the finite element solutions well in the entire parametric space. 

\begin{table}[!htbp]
	\begin{center}
			\centerline{\begin{tabular}{| c | c | c | c | c |}
			\hline
			$   $ & $  $ & $  $ & $  $ & $  $ \Tstrut \\
			$   $ & $  \xi = 10^{-5} $ & $  \xi = 10^{-4} $ & $  \xi = 10^{-3}  $ & $  \xi = 10^{-2}  $ \Bstrut \\ 
			\hline
			$   \chi = 10^{-3}   $ & $  $ & $  $ & $  $ & $  $ \Tstrut \\ $ \mathrm{FE} $ & $ 0.25 \times 10^5
			$  & $ 0.25 \times 10^4
			$ & $ 0.25 \times 10^3
			$ & $ 0.26 \times 10^2 $ \\ $ \Psi(\chi,\xi) $ &  $ 0.25 \times 10^5 $ &  $ 0.25 \times 10^4 $ &  $ 0.26 \times 10^3 $ &  $ 0.26 \times 10^2 $ \\ $ \Psi^{i}(0,\xi) $ & $ 0.25 \times 10^5 $ & $ 0.25 \times 10^4 $ & $ 0.25 \times 10^3 $ & $ 0.25 \times 10^2 $ \\ $ \Psi^{c}(\chi,\xi) $ & $ 1.08 \times 10^7 $ & $ 0.85 \times 10^7 $ & $ 0.62 \times 10^7 $ & $ ~ 0.39 \times 10^7 $ \Bstrut\\ 
			\hline
			$  \chi = 10^{-2} $ & $  $ & $  $ & $  $ & $  $ \Tstrut \\ $ \mathrm{FE} $ & $ 0.12 \times 10^5 
			$  & $ 0.22 \times 10^4 
			$ & $ 0.25 \times 10^3 
			$ & $ 0.26 \times 10^2 $ \\ $ \Psi(\chi,\xi) $ &  $ 0.12 \times 10^5 $ &  $ 0.22 \times 10^4 $ & $ 0.25 \times 10^3 $ &  $ 0.26 \times 10^2 $ \\ $ \Psi^{i}(0,\xi) $ & $ 0.25 \times 10^5 $ & $ 0.25 \times 10^4 $ & $ 0.25 \times 10^3 $ & $ 0.25 \times 10^2 $ \\ $ \Psi^{c}(\chi,\xi) $ & $ 1.08 \times 10^5 $ & $ 0.85 \times 10^5 $ & $ 0.62 \times 10^5 $ & $ ~ 0.39 \times 10^5 $ \Bstrut\\ 
			\hline
			$  \chi = 10^{-1} $ & $  $ & $  $ & $  $ & $  $ \Tstrut \\ $ \mathrm{FE} $ & $ 0.54 \times 10^3
			$  & $ 0.32 \times 10^3
			$ & $ 0.12 \times 10^3
			$ & $ 0.23 \times 10^2 $ \\ $ \Psi(\chi,\xi) $ & $ 0.55 \times 10^3 $ &  $ 0.32 \times 10^3 $ &  $ 0.12 \times 10^3 $ &  $ 0.23 \times 10^2 $ \\ $ \Psi^{i}(0,\xi) $ & $ 0.25 \times 10^5 $ & $ 0.25 \times 10^4 $ & $ 0.25 \times 10^3 $ & $ 0.25 \times 10^2 $ \\ $ \Psi^{c}(\chi,\xi) $ & $ 1.08 \times 10^3 $ & $ 0.85 \times 10^3 $ & $ 0.62 \times 10^3 $ & $ ~ 0.39 \times 10^3 $ \Bstrut\\
			\hline
			$  \chi = 1 $ & $  $ & $  $ & $  $ & $  $ \Tstrut \\ $ \mathrm{FE} $ & $ 1.02 \times 10^1
			$  & $ 0.80 \times 10^1 $ & $ 0.57 \times 10^1
			$ & $ 0.34 \times 10^1 $ \\ $ \Psi(\chi,\xi) $ & $ 1.03 \times 10^1 $ & $ 0.81 \times 10^1 $ &  $ 0.56 \times 10^1 $ &  $ 0.34 \times 10^1 $ \\ $ \Psi^{i}(0,\xi) $ & $ 0.25 \times 10^5 $ & $ 0.25 \times 10^4 $ & $ 0.25 \times 10^3 $ & $ 0.25 \times 10^2 $ \\ $ \Psi^{c}(\chi,\xi) $ & $ 1.08 \times 10^1 $ & $ 0.85 \times 10^1 $ & $ 0.62 \times 10^1 $ & $ ~ 0.39 \times 10^1 $ \Bstrut\\ 
			\hline
	\end{tabular}}
	\end{center}
	\caption{Dimensionless functions $\Psi$ for finite element solutions (FE), $\Psi(\chi,\xi)$, $\Psi^{i}(0,\xi)$, and $\Psi^{c}(\chi,\xi)$.}
	\label{forces for the spheres}
\end{table}

\section{Two-parameter versus single-parameter asymptotic solutions}

The new asymptotic solutions developed in Sections 2 and 3 allow one to identify conditions under which the single-parameter asymptotic solutions, for $\chi=0$ or $\chi = {\cal O}(1)$, are acceptable. This is not just interesting, but also important in deciding whether one can treat the layer material as incompressible, and thus measure just one material constant. This issue is straightforward to address for flat layers, for which  one can simply compare the normalized apparent moduli (\ref{AY incomp}) and (\ref{AY comp}) versus (\ref{YA}). To this end, let us use (\ref{zeta}) to eliminate $\chi$ in favor of the parameter $\zeta$, which has been already exploited for extracting the solution for $\chi= {\cal O} (1)$ in Section 2.3. Then  (\ref{YA}-\ref{AY comp}) yield  the following ratios:
\begin{equation}\label{ratio incompressible}
\frac{\hat{E}^i}{\hat{E}} = \frac{I_0\left(\frac{1}{\zeta
	}\right)}{8 \zeta ^2
	\left[I_0\left(\frac{1}{\zeta }\right)-2 \zeta 
	I_1\left(\frac{1}{\zeta	}\right)\right]} +{\cal O} \left( \xi^2 \right) 
\end{equation}
and 
\begin{equation}\label{ratio compressible}
\frac{\hat{E}^c}{\hat{E}} = \frac{I_0\left(\frac{1}{\zeta
	}\right)}{I_0\left(\frac{1}
	{\zeta }\right)-2 \zeta 
	I_1\left(\frac{1}{\zeta
	}\right)} +{\cal O} \left( \xi^2 \right)  \: .
\end{equation} \\
These expressions make it clear that $\zeta$ is an appropriate single parameter for identifying the domains of validity for the single-parameter asymptotic solutions. This can be done as follows. For a chosen  error tolerance, one can solve  (\ref{ratio incompressible}) for  $\zeta^i$, such that for all $\zeta > \zeta^i$, the error associated with using $\hat{E}^i$ rather than $\hat{E}$ will be below the tolerance. Similarly, one can solve  (\ref{ratio compressible}) for  $\zeta^c$, such that for all $\zeta < \zeta^c$, the error associated with using $\hat{E}^c$ rather than $\hat{E} $ will be below the tolerance. That is, the single-parameter asymptotic solutions become too inaccurate in the interval  $\zeta^c < \zeta < \zeta^i\:. $
For example, let us choose the error tolerance of 10\%, and solve (\ref{ratio incompressible}) for $ \zeta^i  = 1.3$ and (\ref{ratio compressible}) for $ \zeta^c  = 0.046$. Then, for $\xi = 10^{-2}$ and $\xi = 10^{-3}$,  these bounds translate into the intervals $0.49 < \nu < 0.49999$ and $0.4999 < \nu < 0.4999999$, respectively. Thus for thin flat layers, with rare exceptions, the apparent response should not be treated as incompressible. 
 This point has been well appreciated in the solid mechanics literature, where the approximate solutions accounting for compressibility are common. 
 
 \bigskip
 \noindent
 At this point, we can divide the parametric $\zeta$-axis into three intervals $(0,\zeta^c )$, $(\zeta^c ,\zeta^i )$, and $(\zeta^i , \infty )$, and refer to them as compressible, intermediate, and incompressible regimes, respectively. Further, we refer to the transition points  $\zeta^c $ and $\zeta^i $ as nearly compressible and nearly incompressible, respectively. 
 
\bigskip
\noindent
Identification of regimes and transitions for layers between two spheres has to be approached differently,  as we do not have an explicit expression for the force. To this end we observe that $\zeta$ could have been identified directly from  (\ref{equation for A}), as the dimensionless characteristic length for $R$. By applying this logic to    (\ref{equation for A spheres}) we conclude that, for layers between two spheres, the dimensionless characteristic length should be
\begin{equation}\label{characteristing length spheres}
\bar\zeta := \frac{\sqrt{\xi}}{\chi}\:.
\end{equation}
In the absence of equations similar to (\ref{ratio incompressible}) and (\ref{ratio compressible}), we can assess the usefulness of $\bar\zeta$ from Table~\ref{forces for the spheres}, by comparing $\Psi^i$ and $\Psi^c$ with the finite element solutions. We observe that $\Psi^i (0,\xi )$ is accurate if and only  if $\bar\zeta \geq 1$. Thus we can adopt ${\bar\zeta}^i =1$ as the nearly incompressible transition point. In contrast, the data does not support the notion that there is a well-defined ${\bar\zeta}^c$. Rather it appears that the nearly compressible transition is dictated by $\chi$. This issue will be resolved in the next section.

\section{Alternative asymptotic approach}
\label{Lame_asymp}

In this section, we approach the problem for layers between two spheres by developing asymptotic series approximations for Navier's equations directly, without relying on Love-Galerkin's potential. This approach allows us to identify the transitions and regimes; for mathematical details we refer to \cite{Maz'ya,Movchans1}.

\bigskip
\noindent
Let us begin with rewriting Navier's equations  (\ref{equiulibrium r},\ref{equiulibrium z}) in terms of the scaled coordinates in (\ref{scaled coordinates spheres}):
\begin{eqnarray}
  \frac{1}{\xi^2} \frac{\mu}{\lambda} \frac{\partial^2 u_r}{\partial Z^2} + \frac{1}{\xi^{3/2}} \left( 1+ \frac{\mu}{\lambda}\right) \frac{\partial^2 u_z}{\partial Z \partial R} + \frac{1}{\xi}\left( 1+ \frac{2\mu}{\lambda}\right) \left(\frac{\partial^2 u_r}{\partial R^2}
+\frac{1}{R} \frac{\partial u_r}{\partial R} -\frac{u_r}{R^2}\right) =0,   \nonumber \\
 \frac{1}{\xi^2}\left( 1+ \frac{2\mu}{\lambda}\right) \frac{\partial^2 u_z}{\partial Z^2} +  \frac{1}{\xi^{3/2}}\left( 1+ \frac{\mu}{\lambda}\right)  \left(\frac{1}{R} \frac{\partial u_r}{\partial Z} + \frac{\partial^2 u_r}{\partial Z \partial R}\right) + \frac{1}{\xi}   \frac{\mu}{\lambda}
\left(\frac{\partial^2 u_z}{\partial R^2}
+\frac{1}{R} \frac{\partial u_z}{\partial R}\right) =0. \label{scaled Navier}
\end{eqnarray}
Here we use the ratio $\mu/\lambda$ rather than more conventional $ \lambda / \mu$ because we are interested in  $\nu \approx  1/2$, where  $\mu/\lambda \ll 1$.
Accordingly, $\mu/\lambda = {\cal O}(1)$ in the compressible regime, and can be chosen as ${\cal O}(\xi^q )$ with $q>0$ in the nearly incompressible regime. 

\subsection{Compressible regime}
For $\mu/\lambda = O(1)$, the system of equations (\ref{scaled Navier}) can be written in the matrix-operator form, which highlights its asymptotic structure
\begin{equation}\label{operator equation compressible}
\left\{\frac{1}{\xi^2} {\bf L}_0  + \frac{1}{\xi^{3/2}}  {\bf L}_1 +  \frac{1}{\xi}{\bf L}_2
\right\}
{\bf u } = 0.
\end{equation}
Here
\begin{eqnarray*}
{\bf L}_0 &: =& \left(\begin{matrix} 
\frac{\mu}{\lambda}  & 0 \cr  0 & 1+ \frac{2\mu}{\lambda} 
\end{matrix} \right) \frac{\partial^2}{\partial Z^2}\:,\\
{\bf L}_1 &: =& \left(1+\frac{\mu}{\lambda}\right)\left(\begin{matrix} 
	0  & \frac{\partial^2}{\partial Z \partial R} \cr  \frac{\partial^2}{\partial Z \partial R} + \frac{1}{R} \frac{\partial}{\partial Z} & 0
\end{matrix} \right)\:,\\
{\bf L}_2 &: =&  \left(\begin{matrix} 
	\left(1+\frac{2\mu}{\lambda}\right) \left(\frac{\partial^2}{\partial R^2}
	+\frac{1}{R} \frac{\partial}{\partial R} -\frac{1}{R^2}\right) & 0 \cr  0 &  \frac{\mu}{\lambda}
	\left(\frac{\partial^2 }{\partial R^2}
	+\frac{1}{R} \frac{\partial }{\partial R}\right)
\end{matrix} \right)\:,
\end{eqnarray*}
and
$$
{\bf u}= \left(\begin{matrix} 
u_r \cr  u_z
\end{matrix} \right)\:.
$$

\bigskip
\noindent
For   (\ref{operator equation compressible})   the appropriate asymptotic approximation for the displacement field column-vector is
\begin{equation}\label{asymptotic series}
{\bf u} \simeq {\bf u}^{(0)}(R, Z) + \sqrt{\xi} {\bf u}^{(1)}(R, Z) + \xi{\bf u}^{(2)}(R, Z). 
\end{equation}
Now the asymptotic approximation can be constructed as  a recurrence  of boundary-value problems on the scaled cross section
\begin{equation}
{\bf L}_0 {\bf u}^{(n)} = -{\bf L}_1 {\bf u}^{(n-1)}-{\bf L}_2 {\bf u}^{(n-2)}, ~~\mbox{when}~~ |Z| < 1+\frac{1}{2} R^2,
\label{iterations}
\end{equation}
with the boundary conditions
\begin{equation}
{\bf u}^{(n)} = \pm U {\bf e}_z \delta_{n0} ~~ \mbox{as}~~ Z = \pm\left(1+\frac{1}{2} R^2\right).
\label{DBC}
\end{equation}
In (\ref{iterations}) the terms with negative superscript indices must be set equal to zero. 

\bigskip
\noindent
By solving the sequence of problems defined by (\ref{iterations}) and (\ref{DBC}) for $n=0, 1, 2$, we obtain the asymptotic approximation of (\ref{asymptotic series})
in the form
\begin{equation}\label{displacements Sasha compressible}
{\bf u} \simeq \left(\begin{matrix} 
\sqrt{\xi} \frac{\lambda + \mu}{2 \mu}   R \left[\frac{4Z^2}{(2+R^2)^2}-1\right]\\ 
\frac{2Z}{2+ R^2} -\frac{\xi }{3  } \frac{ \lambda}{  \mu} \frac{2-  R^2}{2+ R^2}  \left[ \frac{2Z^2}{(2+R^2)}-1      \right] Z 
\end{matrix} \right) U \:.
\end{equation} \\
This displacement field coincides with that derived in \cite{PTK} for the $\chi={\cal O}(1)$.

\subsection{Nearly compressible transition}
\label{S4_2}

Now we start considering regimes for which the ratio $\mu/ \lambda \ll 1$. Based on elementary analysis of the coefficients in (\ref{scaled Navier}) we conclude that, as $\mu/ \lambda \rightarrow 0 $, the asymptotic structure 
changes qualitatively from that in (\ref{operator equation compressible}) when 
$$
\frac{\mu}{\lambda} = \sqrt{\xi}.
$$
Then the corresponding matrix-operator equation is
\begin{equation}\label{operator equation nearly compressible}
\left\{\frac{1}{\xi^2} {\bf L}_0  + \frac{1}{\xi^{3/2}}  {\bf L}_1 +  \frac{1}{\xi}{\bf L}_2 +\frac{1}{\xi^{1/2}} {\bf L}_3
\right\}
{\bf u } = 0.
\end{equation} \\
Here
\begin{eqnarray*}
	{\bf L}_0 &: =& \left(  \begin{matrix} 
			0  & 0 \cr  0 & 1 
	\end{matrix} \right) \frac{\partial^2}{\partial Z^2}\:,\\
	{\bf L}_1 &: =& \left(\begin{matrix} 
 \frac{\partial^2}{\partial Z^2}  & \frac{\partial^2}{\partial Z \partial R} \cr  
	\frac{\partial^2}{\partial Z \partial R} + \frac{1}{R} \frac{\partial}{\partial Z} & 2  \frac{\partial^2}{\partial Z^2}
	\end{matrix} \right)\:,\\
	{\bf L}_2 &: =&  \left(\begin{matrix} 
		\frac{\partial^2}{\partial R^2}
		+\frac{1}{R} \frac{\partial}{\partial R} -\frac{1}{R^2} &   \frac{\partial^2}{\partial R \partial Z} \cr    \frac{\partial^2}{\partial Z \partial R} + \frac{1}{R} \frac{\partial}{\partial Z}  &  0
	\end{matrix} \right)\:,\\
	{\bf L}_3 &: =&  \left(\begin{matrix} 2 \left(\frac{\partial^2}{\partial R^2}
		+\frac{1}{R} \frac{\partial}{\partial R} -\frac{1}{R^2}\right) & 0 \cr  0 &  
		\frac{\partial^2 }{\partial R^2}
		+\frac{1}{R} \frac{\partial }{\partial R} 
\end{matrix} \right)\:.
\end{eqnarray*} \\
By extending the setting  of (\ref{asymptotic series}--\ref{DBC}) to include ${\bf L}_3$, we obtain
\begin{equation}\label{displacements Sasha nearly compressible}
{\bf u} \simeq \left(\begin{matrix} 
\frac{1}{2}  R \left[\frac{4Z^2}{(2+R^2)^2}-1\right]+\frac{\sqrt{\xi}}{2}  R \left[\frac{4Z^2}{(2+R^2)^2}-1\right]\\ 
\frac{2Z}{2+ R^2} -\frac{\xi }{3  }  \frac{2-  R^2}{2+ R^2}  \left[ \frac{2Z^2}{(2+R^2)}-1      \right] Z 
\end{matrix} \right) U \:.
\end{equation}
This solution is new. Notice that in contrast to (\ref{displacements Sasha compressible}), where $u_r \ll u_z \simeq U $, the displacement field in (\ref{displacements Sasha nearly compressible}) is characterized by $u_r \simeq  u_z \simeq U$. 

\subsection{Nearly incompressible transition}
\label{S4_3}

As $\mu/ \lambda \rightarrow 0 $,   (\ref{scaled Navier}) changes its asymptotic structure 
  qualitatively again when 
$$
\frac{\mu}{\lambda} =  {\xi}.
$$
Then the corresponding matrix-operator equation is
\begin{equation}\label{operator equation nearly incompressible}
\left\{\frac{1}{\xi^2} {\bf L}_0  + \frac{1}{\xi^{3/2}}  {\bf L}_1 +  \frac{1}{\xi}{\bf L}_2 +\frac{1}{\xi^{1/2}} {\bf L}_3 +{\bf L}_4
\right\}
{\bf u } = 0.
\end{equation}
Here
\begin{eqnarray*}
	{\bf L}_0 &: =& \left(  \begin{matrix} 
		0  & 0 \cr  0 & 1 
	\end{matrix} \right) \frac{\partial^2}{\partial Z^2}\:,\\
	{\bf L}_1 &: =&\left(\begin{matrix} 
		0  & \frac{\partial^2}{\partial Z \partial R} \cr  \frac{\partial^2}{\partial Z \partial R} + \frac{1}{R} \frac{\partial}{\partial Z} & 0
	\end{matrix} \right)\:,\\
	{\bf L}_2 &: =&  \left(\begin{matrix} 
		 \frac{\partial^2}{\partial R^2}
		+\frac{1}{R} \frac{\partial}{\partial R} -\frac{1}{R^2} + \frac{\partial^2}{\partial Z^2}  & 0  \cr 0 & 2  \frac{\partial^2}{\partial Z^2}
	\end{matrix} \right) \:,\\
	{\bf L}_3 &: =&  \left(\begin{matrix} 
		0  & \frac{\partial^2}{\partial Z \partial R} \cr  \frac{\partial^2}{\partial Z \partial R} + \frac{1}{R} \frac{\partial}{\partial Z} & 0
	\end{matrix} \right) \:,\\
	{\bf L}_4 &: =&  \left(\begin{matrix} 2 \left(\frac{\partial^2}{\partial R^2}
		+\frac{1}{R} \frac{\partial}{\partial R} -\frac{1}{R^2}\right) & 0 \cr  0 &  
		\frac{\partial^2 }{\partial R^2}
		+\frac{1}{R} \frac{\partial }{\partial R} 
	\end{matrix} \right)\:.
\end{eqnarray*}

\bigskip
\noindent
The procedure, which worked for the two previous cases, does not work here, because neither  $u_r^{(0)}$ nor $u_z^{(0)}$ depend on $Z$, so that the boundary-value problem in (\ref{iterations}) and (\ref{DBC}) becomes degenerate. Further, it becomes apparent that the asymptotic series (\ref{asymptotic series}) must be replaced with
\begin{equation}\label{asymptotic series nearly incompressible}
{\bf u} \simeq \left(\begin{matrix} 
\frac{1}{\sqrt{\xi}} u_r^{(0)}(R,Z) + u_r^{(1)}(R,Z)\\ 
u_z^{(0)}(R,Z)
\end{matrix} \right) U \:.
\end{equation}
With this ansatz, (\ref{operator equation nearly incompressible}) yields the following system of equations for $u_r^{(0)}(R,Z)$ and  $u_z^{(0)}(R,Z)$
\begin{eqnarray}
\frac{\partial^2 u_r^{(0)}}{\partial Z^2}+ \frac{\partial}{\partial R } \left(  \frac{\partial u_r^{(0)}}{\partial R} + \frac{ u_r^{(0)}}{R} +  \frac{\partial u_z^{(0)}}{\partial Z}  \right) =0 \:,  \nonumber\\
\frac{\partial}{\partial Z }  \left(  \frac{\partial u_r^{(0)}}{\partial R} + \frac{ u_r^{(0)}}{R} +  \frac{\partial u_z^{(0)}}{\partial Z}  \right)= 0 \:. \label{system for nearly compressible} 
\end{eqnarray} \\
With the introduction of
\begin{equation}\label{THETA}
\Theta := \frac{\partial u_z^{(0)}}{\partial Z} +  \frac{\partial u_r^{(0)}}{\partial R} + \frac{ u_r^{(0)}}{R} \:,
\end{equation}
we observe that the second equation in (\ref{system for nearly compressible}) implies that $\Theta$ is independent of $Z$, and therefore the first equation can be easily solved:
\begin{equation}\label{Ur0RZ}
u_r^{(0)}(R,Z) = -\frac{\Theta'(R)}{2 } \left[ Z^2 - \left(1+\frac{1}{2}R^2 \right)^2 \right] \:.
\end{equation}
Once $u_r^{(0)}(R,Z)$ is substituted back in (\ref{THETA}), and the boundary conditions 
$$u_z^{(0)} (R,Z)_{  |  Z=\pm \left(1+ \frac{1}{2}R^2 \right)        }  = \pm U $$
are imposed, 
 one obtains the ordinary differential equation 
\begin{equation}
\label{equation for Theta spheres}
\Theta''(R)+ \frac{7 R^2+2}{R^3+2 R} \Theta'(R) -\frac{12}{\left(R^2+2\right)^
	2}  \Theta(R)=-\frac{24
	U}{\left(R^2+2\right)^3}\:.
\end{equation} 
This equation appears to be very similar to (\ref{equation for A spheres}). Indeed, for the chosen $\mu = \lambda \xi$, the corresponding
$$ \frac{\chi^2}{\xi} \simeq 3  \quad{\rm and}\quad \beta= \frac{{\rm i}}{\sqrt{2}}\:,$$
which means that the left-hand sides of (\ref{equation for Theta spheres}) and (\ref{equation for A spheres}) coincide and 
\begin{equation}\label{Theta A}
\Theta(R) = 6 A(R)U \quad {\rm for }\:\: \beta \simeq  \frac{{\rm i}}{\sqrt{2}}\:.
\end{equation}
We could have arrived to this relationship by recognizing that 
$
\Theta = \xi   {\rm div}\bfu \:,
$
and calculating ${\rm div}\bfu$ from (\ref{Love-Galerkin Potential spheres}).

\bigskip
\noindent 
Let us recognize that (\ref{asymptotic series nearly incompressible}) is characterized by the asymptotic structure $\xi$ identical to that developed in \cite{Jeffrey} for the incompressible case. This means that there is no need to consider ratios $\mu / \lambda $ of order less than ${\cal O}(\xi)$, as those ratios would result in the asymptotic structure of the incompressible case.

\subsection{Implications for the transitions}

The two asymptotic solutions, for $\mu/\lambda = \sqrt{\xi}   $ and $\mu/\lambda =  {\xi} $, lend themselves to a better understanding of the transitions introduced in Section 4. Indeed, it is natural to adopt the former as the point of nearly compressible transition, and   the latter as the point of nearly incompressible transition. 
Further,  $\mu/\lambda =  \sqrt{\xi}  $ and $\mu/\lambda =    {\xi} $ imply 
 \begin{equation}\label{zeta bar prime}
\bar\zeta \simeq \frac{1}{\sqrt{3}} \xi^{1/4} \quad {\rm and}\quad \bar{\zeta}= \frac{1}{\sqrt{3}}\:,
\end{equation}
respectively.
The second condition is consistent with the reasoning in Section 4 for using $\bar{\zeta}$ as the parameter for describing the nearly incompressible transition. In particular, it explains why the function $\Psi^i (0,\xi)$ was a good approximation for $\chi=10^{-2}$ for all $\xi$ except for $\xi=10^{-5} $. Indeed, for this case 
$$
\bar{\zeta}=\frac{1}{\sqrt{10}}
$$ \\
is simply insufficiently large. 
The first condition in (\ref{zeta bar prime}) suggests that the nearly compressible transition  should be characterized by
\begin{equation}\label{zeta tilde}
\tilde\zeta := \frac{\xi^{1/4} }{\chi}\:.
\end{equation}
Indeed, this choice is supported by the data in Table~\ref{forces for the spheres}. The single-parameter asymptotic solution for $\chi={\cal O} (1)$ holds if and only if $\tilde\zeta < 1/\sqrt{10}$, that is, for  the last row.

\section{ Discussion}

The asymptotic solutions developed in Sections 2 and 3 allowed us to address several important open issues pertaining to axisymmetric deformation of thin layers constrained by either two rigid plates or two rigid spheres. While our solution for thin layers between plates has the advantage of relying only on the layer thinness and Saint-Venant's principle as the only assumptions, it does not impact significantly the quantitative aspect of the problem. Indeed, for all practical purposes, the problem was quantified in \cite{LS1963++} back in 1963. Nevertheless our analysis introduced the parameter $\zeta$, which represents the interplay between the thinness and compressibility; see (\ref{zeta}). This single parameter is necessary and sufficient for differentiating among the compressible, intermediate, and incompressible regimes. For thin layers between two spheres, the impact was more significant, as we bridged the chasm between the single-parameter asymptotic solutions for compressible \cite{PTK} and incompressible \cite{Jeffrey} layers. Further, we demonstrated that the interplay between the thinness and compressibility requires two parameters $\bar\zeta$ and $\tilde\zeta$; see (\ref{characteristing length spheres}) and (\ref{zeta tilde}). The former allowed us to identify the transition from the nearly incompressible to incompressible regimes, and the latter from the compressible to nearly compressible regimes. These transitions were identified using the asymptotic analysis in Section 5.
 
\bigskip
\noindent
The two problems considered in this work belong to a well-established class of singularly perturbed elliptic boundary-value problems defined on slender domains. A general theory for these problem has been developed in \cite{Maz'ya}, and numerous applications of this theory to boundary-value problems arising in classical elasticity, 
conductivity, and electromagnetism can be found in \cite{Movchans1}.	An important feature of singularly perturbed elliptic boundary-value problems is the presence of boundary layers. In particular, Saint-Venant's principle  for elastic beams, dating back to 1855, can be formally related to exponentially decaying boundary layers. For plates and shells, formal asymptotic analysis gives rise to a combination of exponentially decaying boundary layers associated with local plane strain and  anti-plane shear boundary-value problems defined on scaled transverse sections. In this context, boundary-value problems for thin plates were analyzed in \cite{gregory1, gregory2}, where Saint-Venant's principle was framed in terms of an  exponentially decaying boundary layer in an infinite stripe. The influence of boundary conditions for plates under bending on the exponential decay was examined in \cite{Fried_Dress}. 
	
\bigskip
\noindent	
Asymptotic solutions presented in Sections 2 and 3 constitute leading order approximations, whereas the analysis of Section 5 is based on developing asymptotic series. For Dirichlet boundary conditions, considered in this work, the algorithm for constructing asymptotic series is straightforward,  as the construction  is reduced to a recurrence of boundary-value problems on the scaled cross-section. The problem becomes more technical if Neumann boundary conditions are prescribed on the upper and lower boundaries of the layer. In general, the asymptotic series ansatz 
can be developed using Jordan chains \cite{Maz'ya}. In particular, a Jordan chain  of length   two is required for each in-plane displacement, and a Jordan chain  of length   four is required for the transverse displacement. Among more problem-specific results, we single out Ling's work \cite{ling}, who examined boundary layers in both compressible and incompressible elastic materials. There it was shown that the boundary layer (i) exhibited an exponential decay, (ii) had an insignificant contribution globally, and (iii) mildly depended on Poisson's ratio. These conclusions are in complete agreement with our analysis.

\bigskip
\noindent
Results and techniques presented in this paper  can be extended in several directions. We describe those extensions by dividing them into three groups, in the ascending order of difficulty. 

\bigskip
\noindent
The problem for two spheres can be solved for non-equal spheres, as it was done in both \cite{PTK} and \cite{Jeffrey}. We decided to solve the simpler problem because the formulas were already too long and not particularly insightful. This generalization results in the surfaces $\partial \Omega^\pm$ asymptotically characterized as 
\begin{equation}\label{non-equal spheres}
Z^- = -1 -\frac{1}{2}R^2 \quad{\rm and }\quad Z^+ = 1 + \frac{1}{2\alpha}R^2 \:.
\end{equation}
Here it is assumed that $a$ is the radius of the lower sphere and $\alpha a$ is the radius of the upper sphere. This generalization will give rise to $\Phi$ being a general cubic polynomial in $Z$, and a slightly different but still solvable in hypergeometric functions equation for $A(R)$. Also note one can choose $\alpha=0$ for a plane $\partial\Omega^+$ or  $\alpha<0$ for a "convergent" layer. Further, since the asymptotic analysis reduces the spheres to parabaloids of revolution locally, it is straightforward to extend our analysis to other surfaces, which locally can be approximated by parabaloids of revolution. Those, for example, include power-log cusps.

\bigskip
\noindent
Problems for axisymmetric domains subjected to non-axisymmetric loading, can be solved using trigonometric series,
$$
\frac{1}{2} a_0 (r,z)+ \sum_{n=1}^\infty a_n (r,z)\cos(n \theta ) + b_n (r,z)\sin(n \theta )\:,
$$
in which Love-Galerkin's potential used in this work can be regarded as the term $a_0 (r,z)$ \cite{Love}. Here $\theta$ is the angle of the cylindrical coordinate system $(r,\theta, z)$. This approach may be useful for analyzing thin layers subjected to shear, bending, and twisting \cite{KK}.

\bigskip
\noindent
The asymptotic analysis of Navier's equations can be extended to full-three dimensional, by simply rewriting the governing equations in scaled coordinates, as  it was done in Section 5. In particular, we can generalize the scaling in (\ref{scaled coordinates spheres}) to Cartesian coordinates chosen so that $x_3$ is in the thickness direction,
$$
X_1 = \frac{x_1}{\sqrt{\xi}a}\:, \quad X_2 = \frac{x_2}{\sqrt{\xi}a}\:, \quad X_3 = \frac{x_3}{{\xi}a}\:.
$$
Here $a$ is the smallest radius of curvature for the confining surfaces  $\partial \Omega^-$ and $\partial \Omega^+$, and the minimum distance between these surfaces is set equal to $\xi a$. 
In these coordinates  the gradient operator takes the form
$$
\nabla = \bfe_1 \frac{\partial}{\partial x_1 }+ \bfe_2 \frac{\partial}{\partial x_2 }+ \bfe_3 \frac{\partial}{\partial x_3 }=
\frac{1}{\sqrt{\xi}a}\left( \bfe_1 \frac{\partial}{\partial X_1 }+ \bfe_2 \frac{\partial}{\partial X_2 }+ \bfe_3 \frac{1}{\sqrt{\xi}} \frac{\partial}{\partial X_3  }\right)\:.
$$
This approach may be very difficult to realize, as the challenging ordinary differential equations will be replaced with even more challenging partial differential equations. Nevertheless it may be useful for identifying transitions from the compressible to intermediate, and intermediate to incompressible regimes for general three-dimensional problems for thin layers. 

\bigskip
\noindent
Let us conclude that the asymptotic analyses in Sections 3 and 5 were nicely connected by the relationship (\ref{Theta A}). On the one hand, it provided an insightful interpretation for $A(R)$, and, on the other hand, it allows us to connect (\ref{equation for A spheres}) and (\ref{equation for Theta spheres}).  This connection, valid in the nearly incompressible regime, echoes that examined in \cite{RAS2}. Indeed, both approaches use Navier's equations as the point of departure, and recognize the importance of $\Theta$ being independent of $Z$, as a cornerstone of the solution. In our case,  Navier's equations are restricted to axisymmetric problems, whereas \cite{RAS2, RAS1} considers three-dimensional problems.

\newpage
\noindent
{\Large\bf Acknowledgment }

\bigskip
\noindent
We are grateful to Ken Liechti, Nanshu Lu, and Shutao Qiao for helpful discussions, and Jonathan Zhang for contributing to finite element computations. We appreciate two insightful reviews which helped us to improve the paper. This work was supported by a grant from the National Science Foundation (CMMI 1663551), and the project ${\rm  MCTool^{21}}$ (6305-1452/1490) co-financed by the European Regional Development Fund through the Operational Program for Competitiveness and Internationalization COMPETE 2020, the North Portugal Regional Operational Program NORTE 2020, and by the Portuguese Foundation for Science and Technology FCT under the UT Austin Portugal Program.

\end{document}